\documentclass[11pt]{article}
\usepackage{epsfig}
\input amssymb.sty
\begin{document}
\sloppy

\font\amsy=msbm10

\newtheorem{defn}{Definition}[section]
\newtheorem{example}[defn]{Example}
\newtheorem{question}[defn]{Question}
\newtheorem{prop}[defn]{Proposition}
\newtheorem{thm}[defn]{Theorem}
\newtheorem{lem}[defn]{Lemma}
\newtheorem{cor}[defn]{Corollary}
\newtheorem{conj}[defn]{Conjecture}
\newtheorem{remark}[defn]{Remark}

\newcommand{\dickebox}{{\vrule height5pt width5pt depth0pt}}
\def\NN{\mathbb{N}}
\def\ZZ{\mathbb{Z}}
\def\QQ{\mathbb{Q}}
\def\CC{\mathbb{C}}

\font\amsy=msbm10
\def\IC{\hbox{\amsy\char'103}}  
\def\sg{\sigma}
\def\al{\alpha}
\def\LL{\Cal L}
\def\BB{\Cal B}
\def\AA{\Cal A}
\def\lm{\lambda}
\def\KK{\Cal K}
\def\dim{\mbox{\rm dim }}
\def\sgn{\mbox{\rm sgn}}
\def\om{\omega}
\def\Om{\Omega}
\def\IC{\hbox{\amsy\char'103}}
\font\diff=msbm10
\def\leneq{\hbox{\diff\char'010}}
\def\geneq{\hbox{\diff\char'011}}
\setlength{\parindent}{0pt}

{\Large \bf 
\begin{center}On a class of algebras defined by partitions
\end{center}}
\bigskip

\centerline{\bf 
A.~Regev\footnote{Partially supported by Minerva Grant No. 8441.}} 
\bigskip
\centerline{Department of Mathematics} 
\centerline{The Weizmann Institute of Science, Rehovot 76100, Israel}
\medskip
\centerline{\it ~~regev@wisdom.weizmann.ac.il}
\vskip 0.8cm

\begin{abstract} 
A class of associative (super) algebras is presented, which naturally
generalize both the symmetric algebra $Sym(V)$ and  the wedge algebra
$\wedge (V)$, where $V$ is a vector-space. These algebras are in a
bijection with those subsets of the set of the partitions which are 
closed under inclusions of partitions.
We study the rate of growth of these algebras, then characterize the
case where these algebras satisfy polynomial identities.
\end{abstract}
\medskip

\centerline{2000 Mathematics Subject Classification: }
\centerline{05A17, 05E10, 16S99, 20C30}
\bigskip

\section{Introduction} \label{int}
Throughout this paper let $F$ denote a field of characteristic zero.
Let $V$ be a vector-space, $\dim V=k$,
$\;T^n(V)=V\otimes\cdots\otimes V$
$n$-times, and $T(V)$ is the tensor algebra:
\[
T(V)=\bigoplus_{n}T^n(V).
\]
Both the Lie group $GL(V)$ (or the Lie algebra $gl(V)$) and the 
symmetric group $S_n$ act naturally
on $T^n(V)$, yielding the isotypic decomposition
\begin{eqnarray}\label{eqnstrip}
T^n(V)
=\bigoplus_{\lm\in H(k,0;n)}W_\lm ,
\end{eqnarray} 
where $H(k,0;n)=\{\lm =(\lm_1,\lm_2,\ldots )\mid \lm_{k+1}=0\}$.
In fact, 
$W_\lm\cong V^\lm _k\otimes S^\lm$, where $S^\lm$ is the Specht module
(with character $\chi ^\lm$),
and $V^\lm _k$ is the corresponding $GL(V)$ (or $gl(V)$)
irreducible module, which
is unique - up to an isomorphism.

Recall the {\it symmetric algebra} $Sym(V)$ and the 
{\it wedge algebra} $\wedge (V)$:\\
$Sym(V)=T(V)/I_C$, where $I_C$ is the two sided ideal in $T(V)$ generated
by the elements $x\otimes y - y\otimes x$,
$I_C=\langle x\otimes y - y\otimes x\mid x,\,y\in V\rangle$;

similarly, $\wedge (V)=T(V)/I_E$ where
$I_E=\langle x\otimes y + y\otimes x\mid x,\,y\in V\rangle$.
We show in~\ref{expl1} that 
\begin{eqnarray}\label{}
I_C=\bigoplus_{\lm;\; (1^2)\subseteq\lm} W_{\lm}, \qquad\mbox{and}\qquad
I_E=\bigoplus_{\lm;\;(2)\subseteq\lm } W_{\lm}.
\end{eqnarray}
This leads to the following construction:
Let $\Bbb Y (n)=Par(n)$ be 
the partitions of $n$ and let
$Par =\Bbb Y=\cup_n \Bbb Y(n)$ be all the partitions ($\Bbb Y$ is the 
so called Young graph).
A subset $\Om\subseteq \Bbb Y$ which is closed under inclusions
of partitions is called a {\bf filter}.
Given a subset $\Om\subseteq \Bbb Y$, define 
\[
I_\Om=\bigoplus_{\lm\in \Om } W_{\lm}
\]
It is shown in Section~\ref{LR} that 
the Littlewood-Richardson-rule (LR-rule) implies that 
$\Om$ is a filter if and only if
$I_\Om=\bigoplus_{\lm\in \Om } W_{\lm}$ is a two-sided ideal in $T(V)$.
In such a case
this yields the quotient algebra $A_\Om=A_\Om(V)=T(V)/I_\Om$. 
These algebras -
as well as the analogue superalgebras - are studied in this paper, 
mostly in the case when the dimension of $V$ is finite. \\
\newline
In the `super' case, $V=V_0\oplus V_1$ with $\dim V_0=k$ and $\dim V_1=\ell$; 
again $T(V)=\bigoplus_{n}T^n(V)$ and $S_n$ has a new sign-permutation
action $*$ on $T^n(V)$; also the general-linear-Lie superalgebra 
$pl(V_0,V_1)=pl(k,\ell)$
acts on $T^n(V)$, see~\cite{bereleregev} for the details. 
This yields the following new isotypic decomposition:
\begin{eqnarray}\label{eqnhook}
T^n(V)=\bigoplus_{\lm\in H(k,\ell;n)}W_\lm,
\end{eqnarray}
where $H(k,\ell;n)=\{\lm =(\lm_1,\lm_2,\ldots )\mid \lm_{k+1}\le\ell\}$,
see Section 3 of~\cite{bereleregev}.
Note the difference between $W_\lm$ in
Equations~(\ref{eqnstrip}) and~(\ref{eqnhook}):
in~(\ref{eqnhook})  $W_\lm\cong V^\lm _{k,\ell}\otimes S^\lm$, 
where $V^\lm _{k,\ell}$ is the corresponding $pl(V)$
irreducible module, which
is again unique - up to an isomorphism, see~\cite{bereleregev}.  
The corresponding construction of the algebras is unchanged:
given a filter 
$\Om\subseteq \Bbb Y$,  again let $I_\Om=\bigoplus_{\lm\in \Om } W_{\lm}$;
by Theorem~\ref{thmLRrule}, again $I_\Om$ is
a two sided ideal in $T(V)$,
yielding the quotient algebra $A_\Om=T(V)/I_\Om$. In this paper we study
some of the properties of these algebras, mostly when $\dim V < \infty$. 
%
%
\section{The main results} \label{MR}
In the next section we study the correspondence between 
subsets $\Om\subseteq \Bbb Y$ and subspaces
$I_\Om\subseteq T(V)$. Proposition~\ref{propFiltIde} shows that such
$\Om$ is a filter (i.e. closed under inclusions of partitions) if and only if
$I_\Om$ is a two sided ideal in $T(V)$. In Section~\ref{sec5}
it is proved that filters in $\Bbb Y$ are always finitely generated,
see Theorem~\ref{thm1}.\\
\newline
The filtration $T(V)=\bigoplus_{n}T^n(V)$ induces the 
filtration $A_\Om(V)=A_\Om=\oplus_nA_\Om(n)$. By considering the dimensions
$\;\dim A_\Om(n)\;$ we can talk about the rate of growth of $A_\Om$.
When the dimension of $V$ is finite, it is shown in Section~\ref{secgrowth}
that $A_\Om (V)$ has an exponential growth - which is an integer:
\begin{thm} (See Theorem~\ref{thm5}).
Let $V=V_0\oplus V_1$ be finite dimensional and let $\Om\subseteq \Bbb Y$
be a filter with the corresponding algebra $A_\Om (V)$. Then $A_\Om (V)$
has an exponential rate of growth $\alpha$- which is an integer,
and $0\le \al\le \dim V$.
\end{thm}
When $\al =1$, that rate of growth is polynomial. 
In Sections~\ref{sec6},~\ref{sec7} and~\ref{sec8} we characterizes the 
algebras $A_\Om (V)$ which are p.i., namely 
which satisfy polynomial identities. We prove
\begin{thm}
The algebra $A_\Om (V)$ is p.i.~if and only if it has a polynomial
rate of growth.
\end{thm}
Section~\ref{sec6} treats the `classical' case, namely the case 
$V_1=0$, so $V=V_0$. Theorem~\ref{thm2} gives necessary and
sufficient conditions on $\Om$ - for $A_\Om (V)$ to be p.i., in which case 
$A_\Om (V)$ always satisfies $[x,y]^r=0$ for some $r$. 
\vskip 0.1 truecm
Sections~\ref{sec7} and~\ref{sec8} treat the `super' case. 
Theorem~\ref{thm3} gives necessary and
sufficient conditions on $\Om$ - for $A_\Om$ to be p.i., in which case 
$A_\Om (V)$ always satisfies $h(x)^r=0$ for some $r$. Here $h(x)$ is
any polynomial identity of $E\otimes E$, and $E$ is the infinite dimensional 
Grassmann algebra. In this case it is possible to choose $h(x)=[x,y]^3$ 
hence, again, in the case of p.i.,~$A_\Om(V)$ satisfies a power of the 
commutator $[x,y]$.
\\ 
\newline
We remark that the Littlewood-Richardson rule (LR-rule) is applied, 
in a rather essential way, in the proofs of Proposition~\ref{propFiltIde}
and of Theorems~\ref{thm2} and~\ref{thm3}.\\
\newline
In Section~\ref{sec4} we examine few special cases of such algebras $A_\Om$.
For the classical symmetric algebra $Sym (V)$
we show that $Sym (V) \cong A_\Om (V)$, where $\Om$ is generated by
the partition $(1,1)=(1^2)$.
Similarly for the wedge algebra:  $\wedge (V)\cong A_\Om (V)$
where $\Om$ is generated by the partition $(2)$.
It is well known that $Sym (V)$ is the associated graded algebra
$gr(U(gl(V)))$  of the enveloping algebra $U(gl(V))$
of $gl(V)$. In the `super'-case $V=V_0\oplus V_1$, and
example~\ref{example1} shows that the associated graded algebra
$gr(U(pl(V_0\oplus V_1)))$ of the enveloping algebra of the Lie 
superalgebra $pl(V_0,V_1)$ is also
of the form $gr(U(pl(V_0\oplus V_1)))\cong A_\Om (V)\;$, 
where here $\;A_\Om (V)\cong Sym(V_0)\otimes \wedge (V_1)$ (see for 
example~\cite{scheunert}), and again $\Om$ is generated by the 
partition $(1,1)=(1^2)$.\\
\newline
\section{Filters in $\Bbb Y$  and ideals in $T(V)$} \label{LR}
Let $V=V_0\oplus V_1$, $\dim V_0=k$ and  $\dim V_1=\ell$, and let
$pl(V_0,V_1)=pl(k,\ell)$ denote the corresponding Lie 
superalgebra~\cite{scheunert}.
Notice that the `classical' case is 
obtained by letting
$V_1=0$: $\;pl(k,0)=pl(V_0,0)=gl(V_0)$.

Start with Equation~(\ref{eqnstrip}),
let $\mu\vdash m$ and $\lm\vdash r$ and in $T(V)$ consider
$W_\mu W_\lm\equiv W_\mu\otimes W_\lm$; it is a $pl(k,\ell)$ 
module in a natural way, and we are
interested in the $pl(k,\ell)$-module-decomposition of that module.
The precise decomposition  is given by the LR-rule, a rule which 
arises from the outer multiplication of characters of symmetric groups:

Let $\chi^\mu {\hat \otimes }\chi^\lm$ 
denote the outer-product of the 
characters  $\chi^\mu $ and $\chi^\lm$, then
\[
\chi^\mu {\hat \otimes }\chi^\lm = \sum_{\nu\vdash m+r}c^\nu _{\mu ,\lm }
\chi ^\nu ,
\]
where the coefficients $c^\nu _{\mu ,\lm }$ are given by the LR--rule,
see for example~\cite{sagan}. 
In particular it follows from that rule that if $c^\nu _{\mu ,\lm }\ne 0$
then $\mu ,\lm\subseteq \nu$.

\begin{thm}\label{thmLRrule}

As $pl(k,\ell)$ modules,
\[
W_\mu W_\lm\cong  
(V^\mu _{k,\ell}\otimes V^\lm _{k,\ell})^{\oplus (f^\mu f^\lm)},
\]
and 
\[
V^\mu _{k,\ell}\otimes V^\lm _{k,\ell}\cong  
\bigoplus_{\nu\in H(k,\ell; m+r)}
(V^\nu _{k,\ell})^{\oplus c^\nu _{\mu , \lm} }.
\]
In particular, if $V^\nu _{k,\ell}$
appears in $W_{\lm}W_{\mu}$ then $\lm,\;\mu\subseteq\nu$.
\end{thm}
{\it Proof.} The first claim follows since, as $pl(k,\ell)$ modules,
$W_\mu\cong   ( V^\mu _{k,\ell})^{\oplus {f^\mu}}$ and similarly
for $W_\lm$, see the remark after
Equation~(\ref{eqnstrip}). We prove the second statement. 

Let $M=V^\mu _{k,\ell}\otimes V^\lm _{k,\ell}$, let 
$n=\mid\mu\mid+\mid\lm\mid$ and let $P$ be the matrix
$P=\mbox{diag}(x_1,\ldots ,x_k,y_1,\ldots ,y_{\ell})$.
Then
\begin{eqnarray}\label{eqnTRACE1}
tr_M(P^{\otimes n})=HS_\mu (x_1,\ldots ,x_k,y_1,\ldots ,y_{\ell})
HS_\lm (x_1,\ldots ,x_k,y_1,\ldots ,y_{\ell}).
\end{eqnarray}
By the proof of Theorem 6.30 in~\cite{bereleregev}  
it suffices to show that 
\begin{eqnarray}\label{eqnTRACE2}
tr_M(P^{\otimes n})=
\sum_{\nu\in H(k,\ell ;n)}c^\nu _{\mu , \lm}HS_\nu 
(x_1,\ldots ,x_k,y_1,\ldots ,y_{\ell}).
\end{eqnarray}
By 5.1 of~\cite{remmel} the right-hand-sides of Equations~(\ref{eqnTRACE1})
and~(\ref{eqnTRACE2}) are equal, which completes the proof.\hfill\dickebox  \\
\newline
\begin{defn}\label{defnfilter}
Recall that $\Bbb Y=\cup_n\Bbb Y (n)$ denote the set of all the partitions.
Given a subset $\Om\subseteq \Bbb Y$, define $I_{\Om}\subseteq T(V)$ by
\begin{eqnarray}\label{eqn6}
I_{\Om}=\bigoplus_{\lm\in\Om}W_{\lm}.
\end{eqnarray}
A subset $\Om\subseteq \Bbb Y$ is called {\bf a filter} if it is
closed under inclusions of partitions: if $\mu\in\Om$ and $\mu\subseteq\lm$
then $\lm\in\Om$.\\
Given partitions $\mu^1,\mu^2,\ldots\in\Bbb Y$, let 
$\Om =\langle \mu^1,\mu^2,\ldots \rangle$ denote the filter generated by 
these partitions:
\[
\langle \mu^1,\mu^2,\ldots \rangle =
\{\lm\in \Bbb Y\mid \mu^i\subseteq\lm\quad\mbox{for some $i$}\}.
\]
\end{defn}
\begin{remark}\label{remark3} 
Let $V=V_0\oplus V_1$ with $\dim V_0=k$ and $\dim V_1=\ell$ finite. Let
$\mu$ denote the $k+1\times \ell +1$ rectangle:
$\mu=((\ell +1)^{k+1})$.
By Equation~(\ref{eqnhook}), if $\mu\subseteq\lm$ then $W_\lm =0$.
Given a filter $\Om$, let $\Om_1$ be the filter obtained by adding
all $\mu\subseteq\lm$ to $\Om$, then $I_\Om=I_{\Om_1}$. 
When $\dim V$ is finite, we shall therefore always
assume that that rectangle $\mu$ is in $\Om$.
\end{remark}

A basic and obvious property of such a subspace $I_\Om\subseteq T(V)$ is the
following.
\begin{prop}\label{basic}
If $V^\nu_{k,\ell}$ appears in $I_\Om$ then $W_\nu\subseteq I_\Om$.
\end{prop}
%
The proof of Proposition~\ref{propFiltIde} below requires the following lemma.
\begin{lem}\label{lem3}
Let $\mu$ and $\nu$ be partitions such that $\mu\subseteq\nu$. Then there
exists a partition $\lm$ such that $\chi^\nu$ appears in
$\chi^\mu{\hat\otimes}\chi^\lm$. Moreover, if $\mu,\;\nu\in H(k,\ell )$ 
then also $\lm\in H(k,\ell )$. Here $H(k,\ell )=\cup_n H(k,\ell;n )$.
\end{lem}
{\it Proof}. Let $a_1,\ldots,a_k$ be the lengths of the rows of 
$\nu /\mu$, then $\nu$ appears in
$\chi^\mu {\hat\otimes}\chi^{(a_1)} {\hat\otimes}
\cdots  {\hat\otimes}\chi^{(a_k)}$. Therefore
there is a  $\chi^\lm$ in $\chi^{(a_1)} {\hat\otimes}
\cdots  {\hat\otimes}\chi^{(a_k)}$
such that $\chi^\nu$ appears in
$\chi^\mu{\hat\otimes}\chi^\lm$. 
The second statement follows from the last statement of 
Theorem~\ref{thmLRrule}.\hfill\dickebox  \\
\newline
\begin{prop}\label{propFiltIde}
Let $\Om\subseteq \Bbb Y$ be a subset
with corresponding subspace $I_{\Om}\subseteq T(V)$. Then $\Om$
is a filter if and only if $I_{\Om}$ is a two-sided ideal in $T(V)$.
\end{prop}
{\it Proof.} First, assume $\Om$ 
is a filter and show that $I_\Om$ is an ideal.
It suffices to show the following: Let $\mu\in\Om$ and let $\lm$ be any
partition, then $W_\mu W_\lm\subseteq I_\Om$. This follows from the last
statement of Theorem~\ref{thmLRrule}.
\vskip 0.1 truecm
Next, assume $I_\Om$ is an ideal and show that $\Om$ is a filter.
Let $\mu\in\Omega$, let $\mu\subseteq\nu$ and show $\nu\in\Omega$.
Recall that $T(V)=\oplus_{\lm}W_{\lm}$.
By definition, $W_{\mu}\subseteq I_{\Omega}$. By  
Lemma~\ref{lem3} there exist
a partition $\lm$ such that $\chi^\nu$ appears in 
$\chi^\mu {\hat \otimes }\chi^\lm$ (i.e. $c^\nu _{\mu ,\lm}\ne 0$).
It follows that $V^{\nu}_{k,\ell}$ appears in 
$W_\mu W_\lm \subseteq I_{\Omega} W_\lm \subseteq I_{\Omega}$. 
By Proposition~\ref{basic}
$W_{\nu}\subseteq I_{\Omega}$, and by the definition, $\nu\in\Omega$.
\hfill\dickebox  \\
\newline
\vskip .2 truecm
The algebras $A_\Om (V)$ can now be introduced.
\begin{defn}\label{defA}
Let $V=V_0\oplus V_1$ with the corresponding $W_\lm$'s as in
Equation~(\ref{eqnhook}). Let $\Om\subseteq \Bbb Y$ with the corresponding
subspace $I_\Om = I_\Om (V)$ as in Equation~(\ref{eqn6}). Let $A_\Om (V)$
be the quotient space:
\[
A_\Om (V)=T(V)/I_\Om (V).
\]
Clearly, $A_\Om (V)$ can be identified with the subspace
\begin{eqnarray}\label{eqn16}
A_\Om (V)\equiv \bigoplus _{\lm\not\in\Om}W_\lm .
\end{eqnarray}
\end{defn}

\begin{remark}\label{remark5}
\begin{enumerate}
\item
In particular, the identification~(\ref{eqn16}) implies that if $V\subseteq V'$
then $A_\Om (V)\subseteq A_\Om (V')$.
\item
When $\Om$ is a filter, $I_\Om (V)$ is a two-sided ideal and 
$A_\Om (V)$ is an associative algebra.
\end{enumerate}
\end{remark}
These algebras $A_\Om (V)$ are the subject of this paper.
\section{Some examples} \label{sec4}
Few examples of algebras $A_\Om$ are given below.
We show that in the `classical' case (namely
$V=V_0,\;V_1=0$), both the symmetric
algebra $Sym(V)$ and the wedge algebra $\wedge (V)$ are of the form
$A_\Om (V)$: $\;Sym(V)\cong A_{\langle (1^2)\rangle}(V)$ and
$\wedge (V)\cong A_{\langle (2)\rangle}(V)$.
In the general case, when $V=V_0\oplus V_1$,
$\;A_{\langle (1^2)\rangle}\cong Sym (V_0)\otimes \wedge (V_1)$. Note that
the associated graded algebra
$gr(U(pl(V_0\oplus V_1)))$ of the enveloping algebra of the Lie 
superalgebra $pl(V_0,V_1)$ is also
of that form: $gr(U(pl(V_0\oplus V_1)))\cong  Sym(V_0)\otimes \wedge (V_1)$
\begin{example}\label{expl1}
Let $V=V_0$, $V_1=0$.
\begin{enumerate}
\item
Let $Sym(V)$ be the symmetric algebra of the vector
space $V$. Then $Sym(V)\cong A_{\langle (1^2)\rangle}$,
where $\Om={\langle (1^2)\rangle}=\{\lm\mid\
 (1^2)\subseteq\lm\}$. Thus, as vector spaces, 
$Sym(V)\cong \oplus_nW_{(n)}$.
\item
Let $\wedge (V)$ be the wedge algebra of the vector
space $V$. Then 
$\wedge (V)\cong A_{\langle (2)\rangle}$, 
where $\Om={\langle (2)\rangle}=\{\lm\mid\
 (2)\subseteq\lm\}$. Thus, as vector spaces, 
$\wedge (V)\cong \oplus_nW_{(1^n)}$.
\end{enumerate}
\end{example}
{\it Proof.} We prove part 2. The proof of part 1 is similar. 
Recall that $\wedge (V)=T(V)/I_E$ where
$I_E=\langle x\otimes y + y\otimes x\mid x,\,y\in V\rangle$.
We show that $I_E=\bigoplus_{\lm;\;(2)\subseteq\lm } W_{\lm}$.
Denote $I_E(n)=I_E\cap T^n(V)$. 
$S_n$ acts on $T^n(V)$ - hence on $I_E(n)$ - (from the right, as 
in~\cite{bereleregev}) by permuting
places, and we show first that it maps $I_E(n)$ into itself. 
$I_E(n)$ is spanned by elements of the form 
$x_1\cdots x_{r-1}(x_rx_{r+1}+x_{r+1}x_r)x_{r+2}\cdots x_{r_n}$ (here
$x_1x_2=x_1\otimes x_2$, etc). It suffices to verify for the 
transpositions $(i,i+1)$, and after some obvious reduction, to verify
that $(1,2)$ maps $x_1(x_2v_3+v_3x_2)$ into $I_E$ (and similarly, that 
$(2,3)$ maps $(x_1v_2+v_2x_2)x_3$ into $I_E$). This is clear, since 
\[
(1,2):\;x_1(x_2x_3+x_3x_2)\to x_2x_1x_3+x_3x_1x_2=
\]
\[
(x_2x_1+x_1x_2)x_3 -x_1(x_2x_3+x_3x_2)+(x_1x_3+x_3x_1)x_2
\in I_E.
\]
Let $T_\lm$ be a tableaux of shape $\lm$. To $T_\lm$ corresponds the 
semi-idempotent $e_{T_\lm}=R^+_{T_\lm}C^-_{T_\lm}\in FS_n$.
Here $R_{T_\lm},\,C_{T_\lm}\subseteq S_n$ are the subgroups of 
the $T_\lm$-row and column permutations, with 
\[
R^+_{T_\lm}=\sum_{p\in R_{T_\lm}}p\qquad\mbox{and}\qquad
C^-_{T_\lm}=\sum_{q\in C_{T_\lm}}sgn(q)q .
\]
We apply the following property of the $W_\lm$'s, which is well known:
Let $T_\lm$ be {\it any} tableaux of shape $\lm$, with the corresponding 
semi-idempotent $e_{T_\lm}$, then 
$W_\lm =T^n(V)e_{T_\lm}FS_n$.
\vskip 0.1 truecm
As usual, write $\lm =(\lm_1,\lm_2,\ldots )$.
We verify that $\oplus_{\lm _1 \ge 2}W_\lm\subseteq I_E$. Let $\lm_1\ge 2$
and let $T_\lm$ be a standard tableaux whose first row starts with $1$ and 
$2$. Write  $e_{T_\lm}=R^+_{T_\lm}C^-_{T_\lm}$, then 
$e_{T_\lm}=(1+(1,2))\cdot a$ for some $a\in FS_n$. Given 
${\bar x}=x_1\cdots x_n\in T^n(V)$, we have 
$x_1\cdots x_n(1+(1,2))=(x_1x_2+x_2x_1)x_3\cdots x_n\in I_E(n)$. Since $I_E(n)$
is closed under the $S_n$ action, it follows that 
${\bar x}e_{T_\lm}FS_n={\bar x}(1+(1,2))aFS_n\subseteq I_E(n)$, 
which clearly implies that $W_\lm\subseteq I_E$.
\vskip 0.15 truecm
Next, verify that $I_E\subseteq \oplus_{\lm _1 \ge 2}W_\lm$. Let $\lm =(2)$
with $T_{(2)}$ standard, then $e_{T_{(2)}}=1+(1,2)$. Given 
$x_1x_2\in T^2(V)$, $\;x_1x_2e_{T_{(2)}}=x_1x_2+x_2x_1\in W_{(2)}$.
By the LR-rule it follows that 
\[
x_1\cdots x_{r-1}(x_rx_{r+1}+x_{r+1}x_r)x_{r+2}\cdots x_{r_n}
\in\bigoplus_{\lm;\;\lm_1\ge 2}W_\lm.
\]
Since these are the generators of $I_E$, the above inclusion follows.
This completes the proof of the second example, and the proof of the first
is similar.
\vskip 0.15 truecm
Next we consider the `super' analogues of the previous examples.
\begin{example}\label{example1}
Let $V=V_0\oplus V_1$, let $\Om_1=\langle (1,1) \rangle$ 
be the filter given by the partition $\mu =(1,1)$, and  let
$\Om_2=\langle (2) \rangle$ 
be the filter given by the partition $\nu =(2)$.
We show that 
\[
A_{\langle (1,1) \rangle}(V_0,V_1)
\cong Sym(V_0)\otimes \wedge (V_1),
\]
and similarly 
\[
A_{\langle (2) \rangle}(V_0,V_1)
\cong Sym(V_1)\otimes \wedge (V_0).
\]
It follows from some basic facts in p.i.~theory that both 
$A_{\langle (2) \rangle}(V_0,V_1)$ and $A_{\langle (1^2) \rangle}(V_0,V_1)$ 
satisfy the identity $[[x_1,x_2],x_3]=0$.
\end{example}

{\it Proof.} Fix bases (either finite or infinite)
$t_1,t_2\ldots \in V_0$ and $u_1,u_2\ldots \in V_1$.
By abuse of notation, $t_i,\,u_j\in A_{\langle (1,1) \rangle}$. Also,
$t_i\cdot t_j:=t_i\otimes t_j +I_\Om$ in $T(V)/I_\Om$, and 
similarly for other products. These basis elements satisfy the following three
commutation-relations in $A_{\langle (1,1) \rangle}$:
\begin{enumerate}
\item
$t_it_j=t_jt_i$;
\item
$t_iu_j=u_jt_i$; and
\item
$u_iu_j=-u_ju_i$.
\end{enumerate}
For example we prove 3. Note that here, the $S_n$ action on $T^n(V)$ 
is the same as in~\cite{bereleregev}, and is denoted by $*$. Now,
\[
u_i\otimes u_j+ u_j\otimes u_i=(u_i\otimes u_j)*e_{(1,1)}
\in V^{\otimes 2}*e_{(1,1)}\subseteq
W_{(1,1)}\subseteq I_\Om, 
\]
hence 
$u_iu_j=-u_ju_i$ in $A_\Om$. Here $e_{(1,1)}$ is the semi-idempotent
$e_{(1,1)}=1-(1,2)$. Similarly for the other two relations 1 an 2. 

These commutation-relations imply the isomorphism
$A_{\langle (1,1) \rangle}\cong Sym(V_0)\otimes \wedge (V_1)$.
\vskip 0.25 truecm
Similar arguments show that 
\[
A_{\langle (2) \rangle}(V_0,V_1)
\cong Sym(V_1)\otimes \wedge (V_0).
\]
We consider the classical case (namely $V_1=0$), and
denote $\oplus_nW_{(n)}=A_C$ and $\oplus_nW_{(1^n)}=A_E$. 
As vector spaces, 
$Sym(V)\cong A_C $ and $\wedge (V)\cong A_E$, and these isomorphisms
make $A_C$ and $A_E$ into algebras.
Notice that if the sum 
$\oplus_n$ starts with $n=0$ then these algebras have $1$, while if it
starts with $n=1$, these algebras are without $1$. 
For the next example we introduce the following notation:
\[
A^*_C = \bigoplus_{n\ge 2}W_{(n)}\qquad
\mbox{and}\qquad 
A ^*_E = \bigoplus_{n\ge 2}W_{(1^n)}.
\]
These are ideals in their respective algebras.

\begin{example}\label{}
Let $\Om=\langle (2,1)\rangle$, then we have the following
algebra-isomorphisms:
\begin{eqnarray}\label{eqn8}
A_{\langle (2,1)\rangle}\cong A_C\oplus A^*_E\qquad\mbox{and also}\qquad
A_{\langle (2,1)\rangle}\cong A^*_C\oplus A_E.
\end{eqnarray}
\end{example}
This follows since, as vector spaces, 
\[
A_{\langle (2,1)\rangle}\cong \left (\bigoplus_nW_{(n)}\right )
\oplus \left (\bigoplus_{n\ge 2}W_{(1^n)}\right )=A_C\oplus A^*_E,
\]
and similarly
\[
A_{\langle (2,1)\rangle}\cong \left (\bigoplus_{n\ge 2}W_{(n)}\right )
\oplus \left (\bigoplus_{n}W_{(1^n)}\right )=A^*_C\oplus A_E.
\]
\bigskip

We conclude with few more examples of `classical' algebras $A_\Om$,
namely, $V_1=0,\;V=V_0$.
\begin{example}\label{}
Let $\Om=\langle (a) \rangle\;$ where $\;a >0$. Recall that when $a=2$,
$A_{\langle (2) \rangle}\cong \wedge (V)$.
\begin{enumerate}
\item
$\dim V =k <\infty$.
If $n > k(a-1)$ then $T^n(V)\subseteq I_\Om$. This follows since 
if $\lm=(\lm_1,\ldots,\lm_k)\vdash n > k(a-1)$ then $\lm_1 \ge a$.
It follows that $\dim A_\Om <\infty$. Thus, if
$1\not\in A_\Om$ then $A_\Om$
is nilpotent: $(A_\Om)^{k(a-1)+1}=0$.
\item
$\dim V= \infty$.
Assume $1\not\in A_\Om$. By Remark~\ref{remark5}.1 and by
case 1 above, any finitely generated subalgebra 
is nilpotent. In particular, $A_\Om$ is nil.

\end{enumerate}
\end{example}
\begin{remark}
When $\Om =\langle (2) \rangle\;$, $A_\Om=A_{\langle (2) \rangle}\cong E$,
where $E=E(V)$ is the corresponding Grassmann (Exterior) algebra. 
In particular, $E$ is p.i., satisfying $[[x,y],z]=0$ (even when 
$\dim V=\infty$). Note that $A_{ \langle (1^3) \rangle}$ is not p.i.~since
it does not satisfy the condition $(b^2)\in\Om$ of Theorem~\ref{thm2}.\\
\newline
{\bf Question:} 
\begin{enumerate}
\item
Let $dim V=\infty$ and let $\Om=\langle (3) \rangle\;$.
Is the algebra $A_\Om=A_{\langle (3) \rangle}$ a p.i.~algebra?
\item
Again $\dim V =\infty$, and now $\Om = \langle (2^2) \rangle$. Is
$A_{ \langle (2^2) \rangle}$ a p.i.~algebra? Note that if 
$\dim V < \infty $ then $A_{ \langle (2^2) \rangle}$ is p.i.~by
Theorem~\ref{thm2}.
\end{enumerate}
\end{remark}
%
%


\section{The growth of $\dim (A_\Om(n))$} \label{secgrowth}
It is shown here that for any filter $\Om$, $A_\Om$
has an exponential rate of growth and with an integer exponent. This is
Theorem~\ref{thm5} below.

\begin{defn}
We say that the sequence $d_n$ has exponential rate of growth $\ge \al$
if there exist a polynomial $p$ such that for large enough $n$,
$0< p(n)$ and
$d_n\ge p(n) \al ^n$. In such a case we denote $Exp(d_n)\ge \al$.
Define $Exp(d_n)\le \al$  similarly, and define $Exp(d_n)=\al$ if both 
conditions hold.
Also denote $Exp(d_n)=0$ if $d_n=0$ for all $n$ large enough.
Denote $Exp(A_\Om)=Exp(d_n)$
where $d_n=\dim (A_\Om(n))$. 
\end{defn}
\begin{thm}\label{thm5}
Let $V=V_0\oplus V_1$ be finite dimensional and let $\Om\subseteq \Bbb Y$
be a filter with the corresponding algebra $A_\Om=A_\Om(V)$. Then $A_\Om$
has an exponential rate of growth - which is an integer. More 
precisely, denote $d_n=\dim (A_\Om (n))$, then there exists an
integer $0\le \alpha\le \dim (V_0)+ \dim (V_1)$ such that 
$Exp(A_\Om)=Exp(d_n)=\alpha$.
\end{thm}
{\it Proof.} The proof is given below, see Theorem~\ref{thm6}.

\begin{defn}
\begin{enumerate}
\item
Let $a_1, a_2,b\ge 0$ be integers such that $b\ge a_1, a_2$. Denote by
$D(a_1,a_2,b)$ the following hook-rectangular diagram (i.e.~partition) 
\[
D(a_1,a_2,b)=(b^{a_1},a_2^{b - a_1}).
\]
For example, $(7^3,2^4)=D(3,2,7)$. Note that 
both the arm-length and the leg-length of $D(a_1,a_2,b)$ equal $b$. Also,
$\mid (b^{a_1},a_2^{b - a_1})\mid = (a_1+a_2)b-a_1a_2$.
\item
Let $\Om\subseteq \Bbb Y$ be a filter 
and let $a\ge 0$ be an integer. We say that $\Om$ satisfies the 
$a$-th hook-rectangular condition if there exist a  large enough 
integer $b>0$ such that 
\[
D(a,0,b),\;D(a-1,1,b),\ldots ,\; D(0,a,b)\in \Om.
\]
%
\end{enumerate}
\end{defn}
\begin{remark}\label{remark4}
\begin{enumerate}
\item
If $\;\Om$ satisfies the $a$-th hook-rectangular condition 
(with $b$) then it 
also satisfies the $a+1$-th hook-rectangular condition (with $b+1$). Thus,
for any non-empty filter $\Om$, there exists 
$0\le a$  minimal such that $\Om$ satisfies the $a$-th 
- but not the $a-1$-th - hook-rectangular condition; we denote it
by $a=h_r=h_r(\Om )$. 
\item
By definition, $\Om$ satisfies the $0$-th hook-rectangular condition 
exactly when $\Om = \Bbb Y$, in which case $A_\Om=0$.
\item
$\Om$ satisfies the $1$-th hook-rectangular condition if and only if
$(b),(1^b)\in \Om $ for some $b>0$, hence
if and only if $\dim A_\Om$ is finite, in which case
$A_\Om$ is nilpotent if $1\not\in A_\Om$.
\item
By Theorem~\ref{thm3}, $\Om$ satisfies the $2$-th hook-rectangular condition 
if and only if $A_\Om$ is p.i.
\item
Let $V=V_0\oplus V_1$ with $\dim V_0=k$ and $\dim V_1=\ell$, then by
assumption $\Om$ contains the $k+1 \times \ell +1$ rectangle $\mu^0$:
$\mu^0=((\ell +1)^{k+1})\in\Om$, see Remark~\ref{remark3}. 
We show that $h_r(\Om)\le k+\ell+1$. 
We show in Remark~\ref{remark6} that $h_r(\Om)= k+\ell+1$
exactly when $\Om =\langle \mu^0\rangle$.
Check that $\Om$ satisfies the $k+\ell +1$-th hook-rectangular condition:
If $a_1+a_2=k+\ell+1$ then either $a_1>k$ or  $a_2>\ell$, so for 
large enough $b$, $\mu^0\subseteq D(a_1,a_2,b)$, hence $D(a_1,a_2,b)\in\Om$.
It follows that for any filter $\Om$,  
$\langle \mu^0\rangle\subset \Om\subseteq \Bbb Y$, there exists 
$0\le a\le k+\ell +1$   such that 
$a=h_r=h_r(\Om )$. 
\end{enumerate}
\end{remark}
\begin{lem}\label{lem7}
Given the integers $a_1, a_2\ge 0$, assume $D(a_1,a_2,b)\not \in \Om$ for all 
integers $b\ge a_1, a_2$. Then $Exp(A_\Om)\ge a_1+a_2$.
\end{lem}
{\it Proof.} 	Let $n=\mid (b^{a_1},a_2^{b - a_1})\mid$. Since
$D(a_1,a_2,b)\not \in \Om$, this implies that 
\[
\dim (A_\Om (n))\ge 
\dim W_{ D(a_1,a_2,b)}\ge f^{ D(a_1,a_2,b)}.
\]
The proof now follows from the asymptotic estimates in Section 7 
in~\cite{bereleregev}, which show that for some polynomial $p$,
\[
f^{ D(a_1,a_2,b)}\ge p(n)(a_1+a_2)^n
\]
and $p(n)>0\;$ if $\;n=\mid (b^{a_1},a_2^{b - a_1})\mid\;$ is large enough.
This completes the proof.\hfill\dickebox  \\
\newline
%
The converse is given by the next lemma.
\begin{lem}\label{lem8}
Let $a>0$ be an integer, let  $\Om$ be a filter, and assume  
$\Om$ satisfies the 
$a$-th - hook-rectangular condition:
there exists an integer $b>0$ such that 
\[
D(a,0,b),\; D(a-1,1,b),\ldots ,\;D(0,a,b)\in \Om.
\]
Then $Exp(A_\Om)\le a-1$.
\end{lem}
{\it Proof.} Here is a sketch of the proof. 

Let $\lm\vdash n$. If $\lm\in\Om$ then $W_\lm$ does not contribute
to $\dim A_\Om(n)$. Hence assume $\lm\not\in\Om$. Let $r_1$ (resp. $r_2$)
denote the number of rows (resp. columns) of $\lm$ whose length is
$\ge b$, then $D(r_1,r_2,b)\subseteq \lm$.
If $r_1+r_2 \ge a$ then  by assumption and by Remark~\ref{remark4}.1
$D(r_1,r_2,b)\in\Om$, hence also $\lm\in\Om$, a contradiction.
Since $\lm\not\in \Om$, $r_1+r_2 \le a-1$.\\
By the choice of $r_1$ and $r_2$
it follows that except for its $b\times b$ initial-corner-part, such $\lm$ is
contained in the $(r_1,r_2)$-hook. Since $r_1+r_2 \le a-1$, it follows
that for some polynomial $q$, independent of $\lm$,
$f^\lm\le q(n)(r_1+r_2)^n\le q(n)(a-1)^n$ for all $n$: this follows by a slight
extension of the asymptotic estimates in Section 7 of~\cite{bereleregev}.
Recall that here $W_\lm\cong V^\lm_{k,\ell}\otimes S^\lm$ (see
Equation~(\ref{eqnhook})), hence $\dim W_\lm= \dim (V^\lm_{k,\ell})\cdot
f^\lm$. By~\cite{bereleregev}, $\dim (V^\lm_{k,\ell})$ is polynomialy
bounded, and also the total number of $\lm$'s in the above 
$(r_1,r_2)$-extended-hook, is
polynomialy bounded - as a function of $\mid \lm \mid = n$. 
It therefore follows that for some
polynomial $p$, 
\[
\dim(A_\Om (n)) \le p(n) (a-1)^n.
\]
This completes the proof.\hfill\dickebox  \\
\newline
We can now reformulate and prove Theorem~\ref{thm5}.
\begin{thm}\label{thm6}
Let $V=V_0\oplus V_1$ be finite dimensional, with $\dim V_0=k$ and 
$\dim V_1=\ell$. 
Let $\Om\subseteq \Bbb Y$ be a filter with the corresponding 
algebra $A_\Om$ and let $a=h_r=h_r(\Om)$
as in Remark~\ref{remark4}.1. Assume $\Om\ne \Bbb Y$, 
hence $a\ge 1$. Then $Exp(A_\Om)=a-1$.
\end{thm}
{\it Proof.} As was explained in Remark~\ref{remark4}, such $a=h_r=h_r(\Om)$
exists, and $a\ge 1$ since $\Om\ne 0$. By definition there exist
integers $a_1,a_2\ge 0$ and $a_1+a_2=a-1$, such that for all 
integers $b> a_1,a_2$, $D(a_1,a_2,b)\not\in\Om$. By Lemma~\ref{lem7}
deduce that $Exp(A_\Om)\ge a-1$.\\
Conversely, $\Om$ does satisfy the $a$-th hook-rectangular condition, 
and by Lemma~\ref{lem8} $Exp(A_\Om)\le a-1$, which completes the proof.
\hfill\dickebox  \\
\newline
\begin{remark}\label{remark6}
Clearly, $Exp(T(V))=\dim V_0+\dim V_1=k+\ell$. The converse is also true:
if $Exp(A_\Om )=k+\ell$ then 
$\Om=\langle ((\ell +1)^{k+1})\rangle$ and 
$A_\Om =T(V)$.
Indeed let $\Om\subset \Om_1$, a proper inclusion, and show that 
$Exp(A_{\Om_1})\le k+\ell -1$.
Let $\eta\in\Om_1$, $\;\eta\not \in\Om$. Then $\eta\in H(k,\ell )$.
Now $\eta\subseteq D(k,\ell,\eta _1 +\eta '_1)$, therefore
$D(k,\ell,\eta _1 +\eta '_1)\in\Om_1$. Consider $\lm$'s such that 
$W_\lm$ contribute to $A_{\Om_1}$, namely $\lm\not\in \Om_1$. If
$\lm\not\in \Om_1$, it follows that either $\lm _k\le \eta _1 +\eta '_1$
or $\lm ' _{\ell}\le \eta _1 +\eta '_1$. The asymptotics of $f^\lm$
for such $\lm$'s is $\le (k+\ell-1)^{\mid\lm\mid}$; 
again, this follows by a slight
extension of the asymptotic estimates in Section 7 of~\cite{bereleregev}.
As in the previous 
lemma, this implies that $Exp(A_{\Om_1})\le k+\ell -1$.
\end{remark}
\vskip 0.25 truecm
\section{The `classical' algebras $A_\Om$ which are p.i.} \label{sec6}
Next we characterize those algebras $A_\Om$ that are p.i., 
namely satisfy polynomial identities. Recall that $V=V_0\oplus V_1$. In
this section we consider  the case $V_1=0$.
Recall that $[x,y]=xy-yx$.
We prove:
\begin{thm}\label{thm2}
Let $V_1=0$, so $V=V_0$, and let $\dim V=k<\infty$. 
let $\Om$ be a non-empty filter 
with the corresponding algebra $A_{\Om}=T(V)/I_{\Om}$.
Then the following three conditions are equivalent.
\begin{enumerate}
\item
$A_{\Om}$ is p.i.
\item
$A_{\Om}$ has a polynomial rate of growth.
\item
There exist $c>0$ such that $(c^2)\in\Om$.
In that case $A_{\Om}$ satisfies the identity 
$[x_1,x_2]\cdots [x_{2\ell -1},x_{2\ell}]=0$ (hence also
$[x,y]^{\ell}=0$) where $\ell = c\cdot (\dim V -1 )+1$.
\end{enumerate}
\end{thm}
\medskip
{\it Proof}.
By Theorem 4.13 of~\cite{berele}, 1 implies 2. 
\vskip 0.15 truecm
Show that 2 implies 3: Assume  $(c^2)\not\in\Om$ for all $c>0$
and show that $Exp(A_\Om )\ge 2$, namely, that the growth of $A_\Om$ is larger
than polynomial. Indeed, as in the proof of Lemma~\ref{lem7}, that assumption 
implies that 
\begin{eqnarray}\label{eqn4}
\dim A_\Om(n)\ge \sum_{\lm = (\lm_1,\lm_2)\vdash n}f^\lm \simeq p(n)\cdot 2^n
\end {eqnarray}
for some polynomial $p(x)$, hence $Exp(A_\Om )\ge 2$.
\vskip 0.15 truecm
Finally, show that 3 implies 1 - with the above polynomial identity. 
So, assume that $(c^2)\in \Om$ and show
that $A_{\Om}$ satisfies the identity 
$[x_1,x_2]\cdots [x_{2\ell -1},x_{2\ell}]=0$, where 
$\ell = c\cdot (\dim V -1 )+1$.
 
Let $s_1,s_2\in T(V)$, then 
\begin{eqnarray}\label{eqn1}
[s_1,s_2]\in \bigoplus_{\lm\;;\;|\lm|-\lm_1\ge 1}W_\lm .
\end {eqnarray}
Indeed, fix a basis $v_1,\ldots, v_k\in V$, and without loss of generality
let $s_1=v_{i_1}\cdots v_{i_p}\;(=
v_{i_1}\otimes\cdots\otimes v_{i_p})$ and
$s_2=v_{j_1}\cdots v_{j_q}$.

Now, $[v_i,v_j]\in W_{(1^2)}$ while, for example, 
$[v_{i_1}v_{i_2},v_{j}]=
[v_{i_1},v_{j}]v_{i_2}+v_{i_1}[v_{i_2},v_{j}]\in 
\bigoplus_{\lm;|\lm|-\lm_1\ge 1}W_\lm$, etc. This verifies~(\ref{eqn1}).

By the LR-rule, for $s_1,\ldots,s_{2\ell}\in T(V)$,
\begin{eqnarray}\label{eqn2}
[s_1,s_2]\cdots [s_{2\ell -1},s_{2\ell}]
\in \bigoplus_{\lm\;;\;|\lm|-\lm_1\ge \ell}W_\lm .
\end {eqnarray}
If $\lm$ has more than $k$ parts then by assumption  $\lm\in\Om$ hence
$W_\lm\subseteq I_\Om$.

Assume now that $\lm=(\lm_1,\ldots ,\lm_k)$, with $|\lm |-\lm_1\ge \ell$
and show that $\lm_2>c$. Assume not, then 
since $c\ge\lm_2\ge\ldots\ge\lm_k$ therefore 
$|\lm |-\lm_1 =\lm_2 +\cdots +\lm_k \le (k-1)c<\ell$, a contradiction
(recall: $\ell = c\cdot (\dim V -1 )+1$).
It follows that in that case $\lm_2 > c$ hence, by assumption, 
$\lm\in\Om$, so $W_\lm\subseteq I_\lm$, which is zero in $A_\Om$.
Following equation~(\ref{eqn2}) we see that 
$[s_1,s_2]\cdots [s_{2\ell -1},s_{2\ell}]=0$ in $A_\Om$. This proves
part 2.\hfill\dickebox  \\
\newline
\section{The `super' algebras $A_\Om$ which are p.i.} \label{sec7}
Let $g(x_1,\ldots ,x_r)$ be a multilinear polynomial which is an identity 
of $E\otimes E$, where $E$ is the infinite-dimensional Grassmann
(Exterior) algebra. At the end of Section~\ref{sec8} we discuss
such polynomials
$g(x)$ of low degrees.


\medskip

Here we prove
\begin{thm}\label{thm3}
Let $V=V_0\oplus V_1$ where $\dim V_0,\;\dim V_1 <\infty$, and let $\Om$
be a  filter. 
Then the following three conditions are equivalent.
\begin{enumerate}
\item
$A_{\Om}$ is p.i.
\item
$A_{\Om}$ has a polynomial rate of growth.
\item
$D(2,0,b),D(1,1,b),D(0,2,b)\in\Om$ for some $b>0$. In that case
$A_\Om$ satisfies the identity
\[
q_{tr}(x)=
g(x_1,\ldots ,x_r)g(x_{r+1},\ldots ,x_{2r})\cdots 
g(x_{(t-1)r+1},\ldots ,x_{tr})=0
\]

where $t=b^2$ and $g(x_1,\ldots ,x_r)$ is any multilinear identity 
of $E\otimes E$.
\end{enumerate}
\end{thm}
{\it Proof.} The proof is similar to, but more elaborate than that 
of Theorem~\ref{thm2}.\\
\newline
As in the proof of Theorem~\ref{thm2}, 1 implies 2.\\
\vskip 0.15 truecm
Show that 2 implies 3:
%
%
Assume $(a,1^b)\not\in \Om$ for all positive integers $a,b$. 
Repeating the analogue argument in the proof of Theorem~\ref{thm2}, 
we conclude that 
\begin{eqnarray}\label{eqn4}
\dim A_\Om(n)\ge \sum_{\lm\in H(1,1;n)}f^\lm = 2^{n-1},
\end {eqnarray}
hence
$A_\Om$ is not p.i. ~Similarly if all $(b^2)\not\in\Om$ or if all
$(2^b)\not\in\Om$. This shows that if $A_\Om$ is p.i. then for some
$b>0$, $D(2,0,b),D(1,1,b),D(0,2,b)\in\Om$.
\vskip 0.15 truecm
The proof that 3 implies 1 - and with that identity - is 
similar but more elaborate than the proof of the analogue part 
in Theorem~\ref{thm2}. The proof follows from the following 
claims.
\begin{lem}\label{lem4}
Denote $c(\lm )=\mid\lm\mid-\lm_1=\lm_2+\lm_3 +\ldots$ 

Write $T(V)=\bigoplus_{\lm}W_\lm$. Let $g(x_1,\ldots ,x_r)$ as above and
let $m_1,\ldots ,m_r\in T(V)$, then
\[
g(m_1,\ldots ,m_r)
 \in \bigoplus_{\lm;\;\lm_1,\lm_1'\ge 2}W_\lm =
\bigoplus_{\lm;\;c(\lm ), c(\lm ' )\,\ge 1}W_\lm .
\]
\end{lem}
The proof of this lemma is given in the next section, see 
Corrolary~\ref{cor1}.
\vskip 0.35 truecm
Together with the LR-rule, Lemma~\ref{lem4} implies 
that for any $m_1,\ldots , m_{rt}\in T(V)$,
\[
q_{rt}(m)=g(m_1,\ldots ,m_r)g(m_{r+1},\ldots ,m_{2r})\cdots 
g(x_{(t-1)r+1},\ldots ,x_{tr})
\in\bigoplus_{\lm;\;c(\lm ), c(\lm ' )\,\ge t}W_\lm .
\]
\begin{lem}\label{lem5}
Let $t\ge b^2$ and let $\lm$ be a 
partition satisfying $c(\lm ),c(\lm ' ) \ge t$. Then either 
$(b,1^{b-1})\subseteq \lm $ or $(b^2)\subseteq \lm $ or $(2^b)\subseteq \lm $.
\end{lem}
{\it Proof.} Assume $(b^2),\,(2^b),\,(b,1^{b-1})\not\subseteq \lm $. Since 
$(b^2)\not\subseteq \lm$, $\lm_2<b$. Similarly, $(2^b)\not\subseteq \lm$
implies that $\lm '_2<b$, while $(b,1^{b-1})\not\subseteq \lm $ 
implies that either $\lm_1<b$ or  $\lm '_1<b$,
say $\lm '_1<b$.
%
In such a case, $\lm_{b+1}=0$ and deduce that 
\[
c(\lm )=\lm_2 +\ldots + \lm_b\le (b-1)^2<b^2,
\]
which is a contradiction.\hfill\dickebox  \\
\newline


{\bf The proof that 3 implies 1 (in Theorem~\ref{thm3}).}\\
This last lemma implies that 
\[
\bigoplus_{\lm;\;c(\lm ), c(\lm ' )\,\ge t}W_\lm
\subseteq \bigoplus_{\lm\in\Om}W_\lm ,
\]
which, by Lemma~\ref{lem4} implies that $q_{rt}(x)=0$ in $A_\Om$ 
and therefore is
an identity of $A_\Om$.\\
This completes the proof of Theorem~\ref{thm3}.\hfill\dickebox\\
\newline
%
%
%
%
%
%
\section{A property of the polynomial identities of $E\otimes E$} \label{sec8}
%

Recall that 
\[
e_{(n)}=\sum_{\sg\in S_n}\sg\qquad\mbox{and}
\qquad e_{(1^n)}=\sum_{\sg\in S_n}(-1)^\sg\sg.
\]
In this section we prove
\begin{thm}\label{thmE}
Let $g(x_1,\ldots ,x_d)$ be a multi-linear
polynomial identity of $E\otimes E$, where
$E$ is the infinite-dimensional Grassmann (Exterior) algebra. Let
$V=V_0\oplus V_1$ and let $m_1,\ldots ,m_d $ be monomials in
$T(V)$ such that $m_1\cdots m_d \in T^n(V)$. Also, recall that the 
group algebra $FS_n$ has the super-action $*$ on 
$T^n(V)$~\cite{bereleregev}.
Then
\[
g(m_1,\ldots ,m_d )*e_{(1^n)}=g(m_1,\ldots ,m_d )*e_{(n)}=0.
\]
By Corollary~\ref{cor1} below, this implies that
\[
g(m_1,\ldots ,m_d )\in \bigoplus_{\lm;\;c(\lm ), (\lm ' )\,\ge 1}W_\lm.
\]
\end{thm}
The proof is given below.\\
\newline
{\bf The functions $f_I(\sg)$.}
{\it First, recall the functions $f_I(\sg)$ from Definition 1.1
of~\cite{bereleregev}: let $E=E_0\oplus E_1$ be the usual decomposition
of $E$, let $I\subseteq\{1,\ldots ,d\}$ and let $a_1,\ldots , a_d\in E_0\cup
E_1$ such that $a_i\in E_1$ if and only if $i\in I$; then $f_I(\sg)=\pm 1$
is given by the equation
\begin{eqnarray}\label{eqn9}
a_{\sg (1)}\cdots a_{\sg (d)}=f_I(\sg)a_1\cdots a_d.
\end{eqnarray}
The functions $f_I(\sg)$ also appear naturally in the following context.
Let $m_1,\ldots ,m_d$ be monomials in $x_1,\ldots ,x_n$ such that
$m_1\cdots m_d=x_1\cdots x_n$.
%
Given $\sg\in S_d$, there is a unique permutation $\eta=\eta _\sg\in S_n$ 
such that
$m_{\sg (1)}\cdots m_{\sg (d)}=x_{\eta (1)}\cdots x_{\eta (n)}$.
Let $I=\{1\le i\le d\mid \deg m_i\quad\mbox{is odd}\}$. Then} 
\begin{eqnarray}\label{eqn10}
sgn (\eta )\;=\;sgn (\eta _\sg)\;=\;f_I(\sg). 
\end{eqnarray}
\begin{lem}\label{lem6}
Let 
\[
p(x_1,\ldots ,x_d)=\sum_{\sg\in S_d}\alpha_\sg 
x_{\sg (1)}\cdots x_{\sg (d)}.
\]
Then $p(x_1,\ldots ,x_d)$ is a polynomial identity of $E\otimes E$
if and only if for any pair of subsets $I_1,I_2\subseteq\{1,\ldots ,d\}$,
\[
\sum_{\sg\in S_d}\alpha_\sg f_{I_1}(\sg )f_{I_2}(\sg )=0.
\]
\end{lem}
{\it Proof.} The proof follows straightforward from Equation~(\ref{eqn9}),
since 
\[
a_{\sg (1)}\cdots a_{\sg (d)}\otimes b_{\sg (1)}\cdots b_{\sg (d)}
=f_{I_1}(\sg )f_{I_2}(\sg )(a_1\cdots a_d\otimes b_1\cdots b_d).
\]
{\bf The proof of Theorem~\ref{thmE}.}  
Let $g(x_1,\ldots ,x_d)$ be a multi-linear
polynomial identity of $E\otimes E$ and write
\begin{eqnarray}\label{eqnD}
g(x_1,\ldots ,x_d)=\sum_{\sg\in S_d} \alpha_\sg 
x_{\sg (1)}\cdots x_{\sg (d)}.
\end{eqnarray}
Let $m_1,\ldots ,m_d $ be monomials in
$T(V)$ such that $m_1\cdots m_d \in T^n(V)$,
namely, $m_1\cdots m_d=z_1\cdots z_n$ with $z_1,\ldots ,z_n\in V$.
By embedding $V$ in a large enough vector-space
we may assume w.l.o.g.~that $z_1,\ldots , z_n $ are
linearly independent, as well as super-homogeneous, namely 
$z_1,\ldots ,z_n\in V_0\cup V_1$.
If $m_i=z_{i_1}\cdots z_{i_r}$, then $m_i$ is super-homogeneous of
degree $\delta (m_i)=1$ if the number of $z_{i_j}$'s which are in 
$V_1$ is odd; otherwise $\delta (m_i)=0$.
 It follows that $m_1,\ldots ,m_d$ are also super-homogeneous, and we 
let $I_1$ be the indices $i$ with $m_i$ having super-degree 1:
$I_1=\{1\le i\le d\mid \delta (m_i)=1\}$.
As a monomial, each $m_i$ has a degree, and we let $I_2$ denote
the $i$'s with $m_i$ of odd degree:
$I_2=\{1\le i\le d\mid \deg (m_i)\;\,\mbox{is odd}\}$.
For example, let $x\in V_0$, $y\in V_1$, and let $m=xyyxy$, then 
$\delta (m)=1$ and $\deg (m)=5$.\\
For each $\sg\in S_d$ let  
$\eta=\eta_\sg$ be the unique permutation in $S_n$ such that 
$m_{\sg (1)}\cdots m_{\sg (d)}=z_{\eta (1)}\cdots z_{\eta (n)}$
(see  Equation~(\ref{eqn10})).
It follows from the definition of $f_I(\sg)$ (namely, from
Equation~(\ref{eqn9})) that for $\sg\in S_d$, 
\begin{eqnarray}\label{eqn11}
m_{\sg (1)}\cdots m_{\sg (d)}=f_{I_1}(\sg )\cdot (m_1\cdots m_d)*\eta_\sg
= (m_1\cdots m_d)*(f_{I_1}(\sg )\eta_\sg).
\end{eqnarray}
By Equation~(\ref{eqn10}) 
\begin{eqnarray}\label{eqn11A}
sgn (\eta )\;=\;sgn (\eta _\sg)\;=\;f_{I_2}(\sg).
\end{eqnarray}
By Equation~(\ref{eqnD}) and Equation~(\ref{eqn11})
\begin{eqnarray}\label{eqn12}
g(m_1,\ldots ,m_d)=(m_1\cdots m_d)*
\left ( \sum_{\sg\in S_d} \alpha_\sg f_{I_1}(\sg)\eta_\sg\right ).
\end{eqnarray}
To calculate $g(m_1,\ldots ,m_d)*e_{(1^n)}$, note that 
$e_{(1^n)}=\sum_{\theta\in S_n} sgn(\theta)\theta$. Since $(m*\tau)*\pi=
m*(\tau*\pi)$ (see Lemma 1.5 of~\cite{bereleregev}), we have
\begin{eqnarray}\label{eqn13}
g(m_1,\ldots ,m_d)*e_{(1^n)}=(m_1\cdots m_d)*
\left ( \sum_{\sg\in S_d} \alpha_\sg f_{I_1}(\sg)\eta_\sg
\sum_{\theta\in S_n}sgn(\theta)\theta\right ).
\end{eqnarray}
But by Equation~(\ref{eqn11A})
\begin{eqnarray}\label{eqn14}
\eta_\sg
\sum_{\theta\in S_n}sgn(\theta)\theta
=sgn (\eta_\sg)\sum_{\theta\in S_n}sgn(\theta)\theta
=f_{I_2}(\sg)\sum_{\theta\in S_n}sgn(\theta)\theta,
\end{eqnarray}
hence 
\begin{eqnarray}\label{eqn15}
g(m_1,\ldots ,m_d)*e_{(1^n)}=(m_1\cdots m_d)*
\left ( \sum_{\sg\in S_d} \alpha_\sg f_{I_1}(\sg)f_{I_2}(\sg)
\sum_{\theta\in S_n}sgn(\theta)\theta\right )
\end{eqnarray}
which equals zero by Lemma~\ref{lem6}. This shows that 
$g(m_1,\ldots ,m_d)*e_{(1^n)}=0$. 
The proof that $g(m_1,\ldots ,m_d)*e_{(n)}=0$ is essentially the same,
but with $I_2$ empty, i.e. $f_{I_2}(\sg)=1$ for all $\sg\in S_d$.
\hfill\dickebox  \\
\newline
\begin{cor}\label{cor1}
Let $g(x_1,\ldots ,x_r)$ be a multilinear polynomial which is an identity
of $E \otimes E$, where $E$ is the infinite dimensional Grassmann algebra.
Let $m_1,\ldots ,m_r\in T(V)$ such that $m_1\cdots m_r\in T^n(V)$, then
\[
g(m_1,\ldots ,m_r)
 \in \bigoplus_{\lm\vdash n;\;\lm_1,\lm_1'\ge 2}W_\lm 
\]
\end{cor}

{\it Proof.} We need to show that the components of $g(m_1,\ldots ,m_r)$ 
in both $W_{(n)}$ and in $W_{(1^n)}$ are zero. By~\cite{bereleregev},
$W_{(n)}=(T^n(V))*e_{(n)}$, and similarly 
$W_{(1^n)}=(T^n(V))*e_{(1^n)}$. Thus, $g(m_1,\ldots ,m_r)$ has a component
$\ne 0$ in $W_{(n)}$ if and only if $g(m_1,\ldots ,m_r)*e_{(n)}\ne 0$,
and similarly for $g(m_1,\ldots ,m_r)*e_{(1^n)}$. This implies the proof.
\hfill\dickebox  \\
\newline
%
\begin{remark}\label{remark7}
We conclude this section with two remarks about the identities of $E\otimes E$.
\begin{enumerate}
\item
A. Popov~\cite{popov} showed that $E\otimes E$ satisfies 
the following two identities:
$[\,[x,y]^2,x]=0$ and $[\,[\,[x_1,x_2],[x_3,x_4]\,],x_5]=0$.
Thus, in Theorem~\ref{thmE} we can choose $g(x)$ to be either the 
multilinearization of $[\,[x,y]^2,x]$ (which is of degree 5) or the 
polynomial 
$[\,[\,[x_1,x_2],[x_3,x_4]\,],x_5]$. Moreover, Popov also showed that these 
identities are of minimal degrees (i.e. $E\otimes E$ satisfies no identity 
of degree four), and the above two identities generate all the identities
of $E\otimes E$.
\item
Explicit identities of $E\otimes E$ 
can also be obtained
via cocharacters. The cocharacters of  $E$
are contained in the (1,1) hook,~\cite{olsson}.
Hence, by~\cite{bereleregev2}, 
the cocharacters of  $E\otimes E$ are contained in the (2,2) hook. It
follows that any element of the two sided ideal
$I_{(3,3,3)}\subseteq FS_9$, when realized as a polynomial, is an identity
of $E\otimes E$. This allows the construction of explicit such identities
- of degree 9. 
For example, $s^3 _3[x_1,x_2,x_3]=0$ is such an identity, where
$s_3[x_1,x_2,x_3]$ is the standard polynomial of degree 3. Since 
$s_3[x_1,x_2,1]=[x_1,x_2]$, it follows that $E\otimes E$ satisfies
$[x_1,x_2]^3=0$, and in Theorem~\ref{thmE} we can choose 
$g(x_1,\ldots ,x_d)=g(x_1,\ldots ,x_6)$ to be the multilinearization
of $[x_1,x_2]^3$.
\end{enumerate}
\end{remark}
\vskip 0.2 truecm
\section{Filters in $\Bbb Y$ are finitely generated} \label{sec5}
The following is obvious: let $\Om_1,\;\Om_2\subseteq \Bbb Y$ be two filters,
then  $\Om_1\cup \Om_2$ is a filter, and 
$I_{\Om_1\cup \Om_2}=I_{\Om_1}+I_{\Om_2}$. Similarly for more filters.
For example let $\mu$ be a partition and let
$\Om=\{\lm\in \Bbb Y\mid \mu\subseteq \lm\}$, then $\Om$ is a filter, 
which we denote by  $<\mu>$.
Similarly $<\mu , \;\eta > 
=\{\lm\in \Bbb Y\mid \mu\subseteq \lm \;\mbox{or}\; \eta\subseteq \lm\}$ 
is a filter, and $<\mu , \;\eta >=<\mu >\cup <\eta >$, etc.
\vskip 0.2 truecm
We prove below that every filter $\Om\subseteq \Bbb Y$ 
is finitely generated, i.e.~there exist $r$ and
partitions $\mu^1,\ldots ,\mu^r$ such that 
$\Om=<\mu^1,\ldots ,\mu^r>$. By standard arguments, this is
equivalent to proving a.c.c. (ascending chain condition) on filters.
\begin{defn}\label{def3}
%
%
Recall the notation $H(k,\ell;n)$ for the partitions of $n$
in the $(k,\ell)$ hook:
\[
H(k,\ell;n)=\{\lm =(\lm_1,\lm_2,\ldots )\mid \lm_{k+1}\le \ell\}
\quad\mbox{and}\quad H(k,\ell)=\cup_n H(k,\ell;n).
\]
Thus $\Bbb Y_k=H(k,0)$ is the 
{\it ``strip''} of the partitions with at most $k$ parts. 
Clearly, $\Bbb Y_k$ is the complement of the filter $\Om=<(1^{k+1})>$.
Therefore  a filter that contains  $(1^{k+1})$ is called
 `a filter in $\Bbb Y_k$'.
Similarly, a filter in the hook $H(k,\ell)$ is any filter that 
contains the $(k+1)\times (\ell +1)$ rectangle, i.e~the partition 
$((k+1)^{\ell +1})$.
\end{defn}
\begin{thm}\label{thm1}
Any filter in $\Om\subseteq \Bbb Y$ is finitely generated. 
\end{thm}
{\it Proof}. Let $\mu\in\Om$, then $\mu$ is contained by some 
$(k+1)\times (\ell +1)$ rectangle. It follows that the complement of
$\Om$ is contained by the $(k,\ell)$-hook $H(k,\ell)$. We prove
the theorem under the assumption that $\Om$ is a filter in $\Bbb Y_k $, 
namely the complement of
$\Om$ is contained in the `strip' $\Bbb Y_k =H(k,0)$. The proof of the
general (i.e.~`hook')-case is similar.

The proof of the $k$-strip case is by induction on $k\ge 1$.

The case $k=1$ is obvious: The complement of $\Om$ is contained
in $$\Bbb Y_1=\{(n)\mid n=1,2,\ldots \}.$$ Let $\Om=\Om_1\subset \Om_2$, 
then $(n)\in \Om_2$ for $n$ large enough. Since 
$\Om_2$ is a filter, it follows that its complement is a subset
of the finite set $\{(1),(2),\ldots ,(n-1)\}$, and the proof follows.

Next, assume that the complement of $\Om$ is contained in $\Bbb Y_k$
(namely $(1^{k+1})\in\Om$), and
that the theorem is true for filters whose complements are contained
by the strips $\Bbb Y_r$ when $r<k$.

Assume $\Om=\Om_1\subset \Om_2$, then $\Om_2$ contains a partition $\mu$
with at most $k$ parts: $\mu =(a_1,\ldots ,a_k)$. Denote
$\Om '_2=\langle \mu, (1^{k+1}) \rangle$, so $\Om '_2\subseteq \Om _2$ . 
It is therefore
suffices to prove a.c.c. on chains of filters 
$\Om '_2\subseteq \Om _3 \subseteq \Om _4\subseteq\ldots$
i.e.~that start with $\Om '_2$. Hence we consider the complement of $\Om '_2$.

Let $B=B(\mu)$ denote
the following finite-union of sets of partitions
\[
B(\mu)=\bigcup_{r=0}^{k-1}\left (\bigcup_{0\le b_k\le\ldots\le 
b_{r+1}\le a_{r+1}-1}
\{(x_1,\ldots ,x_r,b_{r+1},\ldots ,b_k)\mid
x_1\ge\ldots \ge x_k\ge b_{r+1}\}\right )
\]
(when $r=0$ the corresponding set
$\{(b_1,\ldots , b_k)\mid 0\le b_k\le\ldots \le b_1\le a_1-1\}$
is finite).
Let $\lm\in \Bbb Y_k$. It is not difficult to see that 
$\lm\not\in B(\mu)\;$ if'f $\;\mu\subseteq\lm$, namely 
$\;a_i\le \lm_i$ for all $i$, and since $\Om '_2$ is a filter, if'f
$\;\lm\in\Om '_2$ . Thus, if $\lm\not\in \Om '_2$ then $\lm\in B(\mu)$, so for
some $0\le r\le k-1$ and some $0\le b_k\le\ldots\le b_{r+1}\le a_{r+1}-1$,
\[
\lm\in\{(x_1,\ldots ,x_r,b_{r+1},\ldots ,b_k)\mid 
x_1\ge\ldots\ge x_r\ge b_{r+1}\}.
\] 

Note that $\{(x_1,\ldots ,x_r,b_{r+1},\ldots ,b_k)\mid 
x_1\ge\ldots\ge x_r\ge b_{r+1}\}$ is isomorphic to $\Bbb Y_r$ under the 
correspondence
\[
(x_1,\ldots ,x_r,b_{r+1},\ldots ,b_k)\leftrightarrow
(x_1-b_{r+1},\ldots ,x_r-b_{r+1}),
\]
and this isomorphism preserves inclusions of partitions. 
Also, here $r<k$.

By induction, each such set satisfies a.c.c (for filters),
hence the above finite union $B(\mu)$ also satisfies that
condition, and the proof of the theorem follows.\hfill\dickebox  \\
\newline

\end{document}